\documentclass[reqno]{amsart}

\usepackage{amsfonts}
\usepackage{bezier}
\usepackage{euscript}
\usepackage{amsfonts,amssymb,amsmath,amsbsy,amsthm,amscd}
\usepackage{xcolor}
\usepackage{graphicx}
\usepackage{color}
\usepackage[T2A]{fontenc}
\usepackage[utf8]{inputenc}
\usepackage[english]{babel}   

\newtheorem{rem}{Remark}
\newtheorem{prop}{Proposition}

\begin{document}

\title[A distributed model of regime-switching diffusions
 ]{
The Kolmogorov forward equation for a distributed model of regime-switching diffusions
}


\author [Bratus]{A.S. Bratus $^1$}
\author[Rozanova]{ O.S. Rozanova$^2$}

\address[1]{Russian university of transport}
\address[2]{Lomonosov Moscow State University}


\email{alexander.bratus@yandex.ru, rozanova@mech.math.msu.su}



\begin{abstract}
For the regime-switching diffusion process with and without advection term we propose an integro-differential equation describing the densities of states continuously distributed over a segment. We demonstrate that there exists a constructive algorithm for  solving the Cauchy problem. We then show that for some initial distributions of states, the solution can be found explicitly. We also discuss how a model with a discrete number of hidden states can be approximated by a model with continuously distributed states.
\end{abstract}

\maketitle

\section{Introduction}
Regime-switching diffusion stochastic processes have a huge number of applications in various fields such as physics, economics, biology, sociology, etc., see, e.g. \cite{Lux} and references therein. Moreover, they have independent theoretical interest and have been studied for a long time, mostly by probability-theoretical methods \cite{Mao_Yuan}, \cite{HSD_book}. For example, in \cite{Yin_Zhu}, \cite{Baran} the authors studied properties of solutions of stochastic differential equations arising from hybrid switching diffusions. The word “hybrid” highlights the coexistence of continuous dynamics and discrete events. The underlying process has two components. One component describes the continuous dynamics, whereas the other is a switching process between a finite number of states, representing discrete events. The generating operator of a random process is associated with deterministic characteristics such as density and  expectation, described by the forward and backward Kolmogorov equations, respectively. These equations are related to each other because they correspond to adjoint operators, but are generally solved differently. The forward equation, written in the divergent form,  is often called the Fokker-Planck equation. The properties of its solutions are key characteristics of a random process \cite{Risken}. They are also the subject of a large modern field of study in partial differential equations and functional analysis, where the requirements for smoothness of the coefficients are relaxed \cite{Bogachev}. In a model containing multiple states, the Kolmogorov equations become a coupled system in which each equation describes the properties of a single state.

 The Fokker-Planck equation for two-state switching diffusion with two states was apparently first written down in a physical context in \cite{Barenblatt} and \cite{two_temp}. In \cite{Aifantis}, various classical problems for such a pair of equations were considered. The most studied models are those containing two states, since an increase in the number of states leads to a more complex system.
This system is linear and can, in principle, be solved using the Fourier transform. However, the inverse Fourier transform can only be found numerically \cite{Lux}.

The main goal of this paper is to demonstrate that by transitioning to a continuum of states, one can obtain a model for which the solution can be obtained explicitly under specific initial conditions and therefore fully analyzed. A finite-state model can be considered as a special case of a continuous model.

In addition to standard switching diffusion, we also consider a model with advection, which can be called the generalized Ornstein-Uhlenbeck model. This model describes mean reversion, and the density function has a nontrivial limit as $t\to \infty$.

 For $M$ states these equations are
\begin{equation*}\label{1}
w_t(t,x,i)=\sum\limits_{j=0}^M \, q_{ij}(x) w(t,x,j)  +b(x,i) w_{x} (t,x,i)+  \frac12 \sigma^2(x,i) w_{xx} (t,x,i), \quad i=0,1,\dots, M,
\end{equation*}
where $Q=q_{ij}(x)$ is Borel measurable and uniformly bounded for all $i, j =0, 1,\dots, M,$ and $x\in\mathbb R$, and satisfies
so called $q$ -property \cite{Yin_Zhang}
\begin{eqnarray*}
&&(i) \quad q_{ij}(x) \ge  0, \quad   x\in{\mathbb R},  \quad  j\ne i;\\
&&(ii) \quad q_{ii}(x)=-\sum\limits_{j\ne i} q_{ij}(x), \quad  x\in{\mathbb R},  \quad i=0,1,\dots,M.
\end{eqnarray*}

The state corresponding $i=0$ is considered as a main state, the others are called hidden states \cite{Lux}.

The respective Kolmogorov forward equations (the Fokker-Plank equations) for the densities $p(t,x,i)$ of the main and hidden states are
\begin{equation}\label{1p}
p_t(t,x,i)=\sum\limits_{j=0}^M \, q^T_{ij}(x) p(t,x,j)  +(b(x,i) w(t,x,i))_{x}+ \frac12 (\sigma^2(x,i) p (t,x,i))_{xx}, \quad i=0,1,\dots, M.
\end{equation}

System \eqref{1p} is considered together with initial data
\begin{equation*}\label{2}
p(t, x, i)|_{t=0} =  \phi(x, i), \quad  x\in{\mathbb R},  \quad i=1,\dots,M.
\end{equation*}

If $\sigma(x,i)=\sigma(x)$, $b(x,i)=b(x)$, then $$P(t,x)=\sum\limits_{i=0}^M p_i$$ is a solution of the heat equation with drift
\begin{equation*}\label{1W}
P_t=- (b(x) P)_x + \frac12 (\sigma^2(x) P)_{xx}.
\end{equation*}

We consider a continuous analog of system \eqref{1p} and assume that the hidden states are distributed continuously over $(0,1]$, so we can consider one integro-differential equation
\begin{equation}\label{3}
p_t=\int\limits_0^1\, K(x,s,\xi) p(t, x,\xi) \, d\xi -(\beta(x,s) p)_{x} +\frac12 (R^2(x,s) p)_{xx},
\end{equation}
where  $p(t,x,s)$ is a function of $x\in\mathbb R$, $s, \in [0,1]$, $t>0$, $R(x,s)$ is a positive function and $K(x,s,\xi)\in L_2(Q)$, $(s, \xi) \in Q=[0,1]\times [0,1]$ for any $x\in \mathbb R$ with the following  analog of $q$ - property:
\begin{eqnarray}
&&(i) \quad K(x,s,s) <  0, \quad    s \in [0,1];\nonumber\\
&& (ii)\quad \int\limits_0^1 K(x,s,\xi)\, ds=0, \quad  \xi \in [0,1].\label{iiK}
\end{eqnarray}
and initial data
\begin{equation}\label{4}
p(t, x, s)|_{t=0} =  \Phi(x, s), \quad  x\in\mathbb R, \quad s \in [0,1].
\end{equation}
Note that the Fichera theory \cite{Fichera} implies that no boundary conditions should be prescribed at $s=0$ and $s=1$.

The kernel $K$ and diffusion coefficient $R$ can be chosen quite arbitrary, nevertheless, if we want to keep the relation with the discrete model, we have to choose
\begin{eqnarray*}
&& K(x,s,\xi) = q_{ji}, \quad  s\in \left[\frac{j}{M},\frac{j+1}{M}\right), \quad  \xi\in \left[\frac{i}{M},\frac{i+1}{M}\right),\quad i,j=0,1,\dots,M,\\
&& \beta(x,s)= b(x,i), \quad  R(x,s)= \sigma(x,i), \quad  s\in \left[\frac{i}{M},\frac{i+1}{M}\right), \quad x\in \mathbb R.
\end{eqnarray*}
For this step-wise $K$, $\beta$ and $R$
\begin{eqnarray*}\label{Kij}
&&q_{ji}= K\left(x,\frac{j+\frac12}{M},\frac{i+\frac12}{M}\right),\quad i,j=0,1,\dots,M,\\
&& b(x,i)=\beta\left(x,\frac{i+\frac12}{M}\right),\quad \sigma(x,i)=R\left(x,\frac{i+\frac12}{M}\right), \quad x\in \mathbb R,\nonumber
\end{eqnarray*}
but we can keep this rule for the continuous kernels and diffusions, serving as an approximation of the step-wise kernels and diffusions by means the Fourier series.



In what follows we consider $ K=K(s,\xi)$, $ \beta=b(s)x+c(s)$ and  $ R=R(s)$.
This choice of $\beta$ at $b<0$ implies mean reversion and makes the process similar to the Ornstein-Uhlenbeck process..

In this case  the respective backward Kolmogorov equation is similar to \eqref{3} and can also be solved by the method described below.

\section{Construction of a solution}
\subsection
{Formal scheme.}

We take the Fourier transform with respect to $x$ and obtain for $\hat p (t,\mu,s)$
\begin{eqnarray*}\label{5}
&&\hat p_t(t,\mu,s)=\\
&&\int\limits_0^1\, K(s,\xi) \hat p(t,\mu,\xi) \, d\xi - \left(\frac12   R^2(s)\mu^2 + i c(s) \mu\right) \hat p(t,\mu,s) +b(s)\mu \hat p_\mu(t,\mu,s).\nonumber
\end{eqnarray*}

1. If $b(s)=0$, then we change
\begin{equation*}\label{V}
V(t,\mu,s)=\hat p(t,\mu,\xi) \, \exp {(\frac12   R^2(s)\mu^2+i c(s) \mu) t}
\end{equation*}
and get a first order equation
\begin{equation}\label{6}
V_t(t,\mu,s)=\int\limits_0^1\, K(s,\xi) V (t,\mu,\xi) \, d\xi.
\end{equation}

2.
If $b(s)\ne 0$, we make a change
\begin{equation*}\label{Vb}
V(t,\mu,s)=\hat p(t,\mu,\xi) \, \exp {\left(-\frac{R^2(s)}{4 b(s)}\,(1-e^{2b(s) t}) \, \mu^2 - i \frac{c(s)}{{b(s)}}\, (1-e^{b(s) t}) \,\mu\right) }
\end{equation*}
and obtain
\begin{equation*}\label{6_1 }
V_t(t,\mu,s)=\int\limits_0^1\, K(s,\xi) \, V (t,\mu,\xi) \, d\xi +b(s) \mu V_{\mu}(t,\mu,s).
\end{equation*}
Then we
introduce a new independent variable $\tilde \mu = \sqrt{|b(s)|}\mu$ and obtain
\begin{equation}\label{6tilde}
V_t(t,\tilde \mu,s)=\int\limits_0^1\, K(s,\xi) \, V (t,\tilde\mu,\xi) \, d\xi +{\rm sign} (b(s)) \tilde\mu V_{\mu}(t,\tilde\mu,s).
\end{equation}

Let us choose in the space $L_2(0,1)$ an orthonormal basis   $X_l(s)$, $l=0, 1,\dots$, such that $X_0\equiv 1$.

Assume
\begin{equation}\label{7}
 K(s,\xi)= \sum\limits_{n,m=0}^\infty A_{nm} \, X_n(s) \, X_m(\xi),
\end{equation}
due to \eqref{iiK} we have $A_{0 m}= 0$.

We look for the solution of \eqref{6} and \eqref{6tilde} as
\begin{equation*}\label{8}
 V(t, \tilde \mu, s)= \sum\limits_{l=0}^\infty a_{l}(t,\tilde \mu) \, X_l(s).
\end{equation*}
Multiply \eqref{6} (or \eqref{6tilde}) by $X_k(s)$ and integrate over $[0,1]$, $k\in \{0\}\cup\mathbb N$.

 Thus,
for $b(s)=0$,
\begin{equation}\label{9}
 (a_{k}(t, \mu))_t= \sum\limits_{l=0}^\infty A_{k l} a_{l}(t, \mu), \,k\in \{0\}\cup \mathbb N.
\end{equation}
and for $b(s)\ne 0$,
\begin{equation}\label{9b}
 (a_{k}(t,\tilde \mu))_t= \sum\limits_{l=0}^\infty A_{k l} a_{l}(t,\tilde \mu)+{\rm sign}(b(s))\tilde \mu (a_{l}(t,\tilde \mu))_{\tilde \mu}, \,k\in \{0\}\cup \mathbb N.
\end{equation}
 Since  $A_{0 l}= 0$, then $a_{0}(t,\mu)=a_{0}(0, \mu)$, as well as $a_{0}(t,\tilde \mu)=a_{0}(0,\tilde \mu)$.

 If  $a_{l}(t, \mu)$  (or, respectively, $a_{l}(t,\tilde \mu)$) are known, then for $b=0$
\begin{equation*}\label{10}
 \hat p (t,\mu,s)= e^{-(\frac12  R^2(s)\mu^2+i c(s) \mu) t}\, \sum\limits_{l=0}^\infty a_{l}(t, \mu) \, X_l(s)
\end{equation*}
and for $b\ne 0$
\begin{equation*}\label{10}
 \hat p (t,\mu,s)= \, e^{\frac{R^2(s)}{4 b(s)}\,(e^{2b(s) t}-1) \, \mu^2 + i \frac{c(s)}{b(s)}\, (e^{b(s) t}-1) \,\mu}\, \sum\limits_{l=0}^\infty a_{l}\left(t, {\sqrt{|b(s)|}{\mu}}\right) \, X_l(s),
\end{equation*}
and the inverse Fourier transform with respect to $\mu$ gives
\begin{equation}\label{10}
 p (t,x,s)\,= \, \sum\limits_{l=0}^\infty B_{l}(t,x,s) \, X_l(s),
\end{equation}
where for $b=0$
\begin{equation}\label{11}
 B_{l}(t,x,s) = \frac{1}{2\pi}\,\int\limits_{\mathbb R} \,e^{-(\frac12   R^2(s)\mu^2+i c(s) \mu) t}\, a_{l}\left(t, {\mu}\right)\, e^{i \mu x}\, d\mu.
\end{equation}
and for $b\ne 0$
\begin{equation}\label{11b}
 B_{l}(t,x,s) =  \frac{1}{2\pi}\,\int\limits_{\mathbb R} \, e^{{\frac{R^2(s)}{4 b(s)}\,(e^{2b(s) t}-1) \, \mu^2 + i \frac{c(s)}{b(s)}\, (e^{b(s) t}-1) \,\mu}}\,  a_{l}\left(t,{\sqrt{|b(s)|}\mu}\right) \, e^{i \mu x}\, d\mu.
\end{equation}

Let the data \eqref{4} be
\begin{equation}\label{14}
 \Phi(x,s)=\sum\limits_{k=0}^\infty\, g_k(x)\, X_l(s), 
\end{equation}
Then
\begin{equation}\label{13}
 \hat \Phi(\mu,s)= \sum\limits_{k=0}^\infty\, \hat g_k(\mu) \, X_l(s), \quad a_{k}(0,\mu)=\hat g_k(\mu). 
\end{equation}

Therefore $a_k(t,\mu)$   can be found from the linear ODE system \eqref{9}  with coefficients depending on $\mu$ as a parameter for $b=0$. For $b\ne 0$ functions $a_k(t,\tilde\mu)$  can be found from the linear first order PDE system \eqref{9b}.

 Thus, we find coefficients  $ B_k(t,x,s)$  from initial data by  \eqref{11} or \eqref{11b}.

If the coefficients $g_k(x)$ are such that $B_k(t,x,s)$ can be calculated analytically, then $p(t,x,s)$ can be represented as a finite sum or series. Below we present such cases, which are interesting both in themselves for demonstrating the properties of the solution and as reference  for testing numerical results.


\subsection{Existence and uniqueness.}

The existence of classical solution to the Cauchy problem   for the  system of the Kolmogorov backward equations with $C^2$ - smooth coefficients is obtained in \cite{Yin_Zhu}. 

With regard to our problem, based on previously known results, we can formulate the following statement.

\begin{prop}\label{P1}
 Let  $b(s)$, $c(s)$,  $R(s)>0$, $K(s,\xi)$ and $\Phi(x,s)\ge 0$ are bounded and  Lipschitz continuous functions on $[0,1]$, $[0,1]\times [0,1]$ and ${\mathbb R}\times [0,1]$, respectively.
 Then there exist a unique classical nonnegative solution to the problem \eqref{3}, \eqref{4}. It  is bounded, continuous on $[0,T] \times {\mathbb R}\times [0,1]$,  uniformly in $t$ Lipschitz continuous with respect to in $(x,s)$ and, for fixed $s$, $C^1$  in $t$ and $C^2$ in $x$.
\end{prop}

\proof
The result is a particular case of Theorem 2.4 \cite{PIDE}, taking into account arguments of \cite{Baran}, \cite{HSD_book} where it is noticed that the evolution of the discrete component can be represented as a stochastic integral with respect to a Poisson random measure. In the new situation, the function $h$ under the stochastic integral sign  in (7), \cite{Baran}, is not  integer-valued, but takes values in $[0,1]$. The nonnegativity follows from the maximum principle, e.g.\cite{MaxPrinc}. Note that Theorem 2.4 \cite{PIDE} is proved for the backward Kolmogorov equation, but due to the independence of the coefficients from $x$, the results can be applied to the forward Kolmogorov equation as well.
$\square$

\medskip

 Proposition \ref{P1} claims that the regularity in $x$ and $s$ are different. This reflects the fact that the diffusion component has a smoothing effect only in the $x$ direction.

 In some examples we deal with below, the kernels $K(s,\xi)$ and initial data $\Phi(x,s)$ turn out to be only piecewise smooth in $s$ in $\xi$ in the whole $[0,1]$. However, we can use the fact that the equation \eqref{3} does not require any boundary condition at $s=s_0$, $s_0={\rm const}$ and consider in the strip ${\mathbb R}\times (s_1,s_2),$ $0\le s_1<s_2\le 1$. Thus we can prove the existence of a solution to the Cauchy problem for every strip, such that $K(s,\xi)$ is smooth in $ (s_1,s_2)\times (s_1,s_2) $ and $\Phi(x,s)$ is smooth in ${\mathbb R}\times (s_1,s_2)$.

\section{One hidden state: $M=1$}\label{2st}
Assume that the discrete matrix is
\begin{equation}\label{qji}
q_{ji}=\left(\begin{array}{cc}-\lambda_1
 & \lambda_2\\ \lambda_1
 & -\lambda_2
  \end{array}\right),\quad \lambda_i>0, \quad i,j=0,1.
\end{equation}

The kernel $K(s,\xi)$, based on \eqref{qji}, is
\begin{eqnarray*}\label{Ksxi}
&& K(s,\xi) =\left\{\begin{array}{rc}-\lambda_1,\, s\in (0,\frac12),\,  \xi\in (0,\frac12)\\
 \lambda_2, \, s\in (0,\frac12),\,  \xi\in (\frac12,1)\\ \lambda_1,\, s\in (\frac12,1),\,  \xi\in (0,\frac12)
 \\ -\lambda_2, \, s\in (\frac12,1),\,  \xi\in (\frac12,1)
  \end{array}.\right.
 \end{eqnarray*}

 \subsection{Approximation by two first modes.}\label{Stwomodes}

 Let $X_n(s)$, $n\in \{0\}\cup \mathbb N$, be an orthonormal basis in $L_2[0,1]$, such that $X_0\equiv 1$.
 The kernel $K(s,\xi)$ can be expanded into series \eqref{7}, the truncation at the first step is
   \begin{equation}\label{K2}
K(s,\xi)=A_{00}+A_{10} X_1(s) +A_{01} X_1(\xi) + A_{11} X_1(s)  X_1(\xi), \quad A_{00}=A_{01}=0. 
\end{equation}


We consider a distributed generalisation of the system with one hidden state, based on \eqref{K2}. 

The problem \eqref{9}, \eqref{13} reduces to
\begin{eqnarray*}\label{15}
 \dot a_0=0, \quad  \dot a_1=A_{10} a_0 + A_{11} a_1, \quad a_0(0)=\hat g_0(\mu), \quad a_1(0)=\hat g_1(\mu),
\end{eqnarray*}
and can be easily solved. For $b=0$ we have
\begin{eqnarray}\label{16}
  a_0=\hat g_0(\mu), \quad  a_1(t,\mu)= e^{A_{11} t} \left(\hat g_1(\mu)+\frac{A_{10} \hat g_0(\mu)}{A_{11}}\right) - \frac{A_{10} \hat g_0(\mu)}{A_{11}},
\end{eqnarray}
for $b\ne 0$
\begin{eqnarray}\label{16b}
 && a_0(\tilde\mu)=\hat g_0\left( \frac{\mu}{\sqrt{|b(s)|}}\right), \\
 &&a_1(t,\tilde\mu)= \left(\int\limits_0^t\,A_{10} \hat g_0\left(\frac{\mu}{\sqrt{|b(s)|}} e^{-\sigma (t-\eta)}\right) e^{-A_{11}\eta }\, d\eta +   \hat g_1\left(\frac{\mu}{\sqrt{|b(s)|}} \, e^{-\sigma t}\right)\right) e^{A_{11}t},\nonumber
\end{eqnarray}
where $\sigma={\rm sign} \,b(s)$.

\medskip

\subsection{Examples of orthonormal systems} \label{Xex}

1. The most common variant for an orthonormal system consists of trigonometric functions $X_0(s)=1$, $X_n(s)=\sqrt{2}\,\cos \pi n s$, $n\in \mathbb N $. In this case
\begin{equation*}\label{A}
A_{00}=A_{01}=0, \quad  A_{10}=\frac{2}{\pi}(\lambda_2-\lambda_1), \quad  A_{11}=-\frac{8}{\pi^2}(\lambda_1+\lambda_2)<0.
\end{equation*}
$A_{10}$ can have any sign or be zero. Note that $\Big|\frac{A_{10}}{A_{11}}\Big|<\frac{\pi}{4}.$

However, in this case, both at the stage of approximating the piecewise constant kernel \eqref{K2},
we lose the accuracy of the calculations, and to improve it, it is necessary to consider   higher-order terms of the expansion.

\medskip

2. A much more convenient option for the case of simulating discrete models with continuous ones is the Haar wavelet basis \cite{Haar}, \cite{Haar_book}.
Recall that it consists of functions $H_0(s)=1$, and functions $\psi_{ij}(s)=2^{\frac{i-1}{2}} \phi \left(2^{\frac{i-1}{2}} s-j \right)$, $i\in\mathbb N$, $j=0,1,\dots, 2^{i-1}$,
$$\phi(s)=\left\{\begin{array}{ll}1,& s\in [0, \frac{1}{2})\\-1,& s\in [ \frac{1}{2}, 1] \\ 0, &  s\notin[0,1]\end{array} \right..
$$
Functions $\psi_{ij}(s)$ can be sequentially ordered as $H_n(s)$, $s\in\mathbb N$ to obtain together with $H_0(s)$ an orthonormal basis in the space
$L_2(0,1)$. Therefore $H_1(s)=\psi_{10}(s)=\phi(s)$.
In this case
\begin{equation*}\label{AA}
A_{00}=A_{01}=0, \quad  A_{10}=\frac{1}{2}(\lambda_2-\lambda_1), \quad  A_{11}=-\frac{1}{2}(\lambda_1+\lambda_2)<0.
\end{equation*}
$A_{10}$ can have any sign or be zero, $\Big|\frac{A_{10}}{A_{11}}\Big|<1.$

\subsection{Examples of solutions}
Below, we consider two types of initial data and study their dynamics under switching diffusion. In the first case, we consider a Gaussian initial distribution common to both the ground and hidden states. In the second case, the initial distribution is also Gaussian for each state, but characterized by different parameters. In particular, this may be the limiting case of a Gaussian distribution, that is, a delta function. We consider both the case with and without the advective term $\beta$.

\subsubsection{Uniform Gaussian initial distribution}\label{S321}

Let the the initial distribution be uniform in $s$ and
\begin{eqnarray}\label{IDGauss}
 \Phi(x,s)= \frac{1} {\sqrt{\pi}}\,e^{-x^2}.
\end{eqnarray}
This function is normalized such that $\int_{\mathbb R}\, \Phi(x,s)\,dx=1$, $s\in[0,1]$.
It is easy to calculate that $\hat g_0(\mu)=  \, e^{-\frac14 \mu^2}$, $\hat g_1(\mu)=0$. Therefore $a_0(\mu)=\hat g_0(\mu)$, $a_1=\frac{A_{10}}{A_{11}}(e^{A_{11}t}-1) a_0(\mu)$.

\medskip

a. For $b=0$ \eqref{10} and  \eqref{16} give
\begin{eqnarray*}\label{17e}
&& p(t,x,s)=\\&&\frac{1}{\sqrt{\pi (1+2 t R^2(s))}}\, \left(1+ \frac{A_{10}}{A_{11}}\, \left(e^{A_{11} t}-1\right)\,X_1(s)\right)\,e^{-\frac{(x-c(s) t)^2}{1+2 R^2(s)t}}.\nonumber
\end{eqnarray*}

\medskip

b. For $b\ne 0$ \eqref{11b} and  \eqref{16b} give

\begin{eqnarray*}\label{17eb}
&& p(t,x,s)=B_0(t,x,s) + B_1(t,x,s)\,X_1(s),
\end{eqnarray*}
\begin{eqnarray*}\label{17B0}
&&B_0(t,x,s)=
\frac{1}{\sqrt{\pi\left(\frac{R^2(s)}{b(s)}(e^{2 b(s) t}-1)+1\right)}}\,\exp\left(-\frac{\left(x+(1-e^{b(s)t})\frac{c(s)}{b(s)}\right)^2}{\frac{R^2(s)}{b(s)}( e^{2 b(s) t-1})+1}\right).
\end{eqnarray*}
The coefficient $B_1(t,x,s)$ cannot be found explicitly but can be estimated.

Indeed,
\begin{eqnarray}\label{B1}
&&  B_1(t,x,s)\,= \,\frac{A_{10}}{\sqrt{\pi}}\,\int\limits_0^t\,e^{A_{11}(t-\eta)}\,J(\eta;t,x,s)\,d\eta,\\
&& J(\eta;t,x,s)= \frac{1}{\sqrt{-\frac{R^2(s)}{b(s)}(1- e^{2 b(s) t})+
e^{-2\sigma (\eta-t)}}}\,\exp\left(-\frac{\left(x+(1-e^{b(s)t})\frac{c(s)}{b(s)}\right)^2}{-\frac{R^2(s)}{b(s)}(1- e^{2 b(s) t})+e^{-2\sigma (\eta-t)}}\right).\nonumber
\end{eqnarray}
Let us denote
\begin{eqnarray*}
&& J_1(t,x,s)= \frac{1}{\sqrt{-\frac{R^2(s)}{b(s)}(1- e^{2 b(s) t})+
1}}\,\exp\left(-\frac{\left(x+(1-e^{b(s)t})\frac{c(s)}{b(s)}\right)^2}{-\frac{R^2(s)}{b(s)}(1- e^{2 b(s) t})+e^{2\sigma t}}\right),\\
&& J_2(t,x,s)= \frac{1}{\sqrt{-\frac{R^2(s)}{b(s)}(1- e^{2 b(s) t})+
e^{2\sigma t}}}\,\exp\left(-\frac{\left(x+(1-e^{b(s)t})\frac{c(s)}{b(s)}\right)^2}{-\frac{R^2(s)}{b(s)}(1- e^{2 b(s) t})+1}\right).
\end{eqnarray*}
Then for $\sigma=-1$
\begin{eqnarray*}
&& J_1(t,x,s)\le J(\eta;t,x,s)\le J_2(t,x,s),
\end{eqnarray*}
and for $\sigma=1$
\begin{eqnarray*}
&& J_2(t,x,s)\le J(\eta;t,x,s)\le J_1(t,x,s).
\end{eqnarray*}

Since the integral \eqref{B1} with the change of $J$ to $J_1$ or $J_2$ can be computed analytically, we get
\begin{eqnarray*}
&& B_1(t,x,s) = \frac{A_{10}}{A_{11}} (e^{A_{11} t}-1)  {Q(t,x,s)}{\sqrt{\pi}},
\end{eqnarray*}
where
\begin{eqnarray*}
&&\frac{1}{\sqrt{\pi}}\,J_1(t,x,s)\le Q(t,x,s)\le\frac{1}{\sqrt{\pi}}\, J_2(t,x,s), \quad \sigma=-1, \\
&&\frac{1}{\sqrt{\pi}}\,J_2(t,x,s)\le Q(t,x,s)\le \frac{1}{\sqrt{\pi}}\,J_1(t,x,s), \quad \sigma=1.
\end{eqnarray*}

Note that if $b(s)<0$ ($\sigma=-1$) there exists a nontrivial steady state, a limit $ p^*(x,s)=\lim\limits_{t\to \infty}\, p(t,x,s)$.
Namely, if $B_0^*(x,s)=\lim\limits_{t\to\infty}\, B_0(t,x,s)$,  $Q^*(x,s)=\lim\limits_{t\to\infty}\, Q(t,x,s)$, then
\begin{eqnarray*}
 B_0^*(x,s)=
\frac{1}{\sqrt{\sqrt{\pi}\left(-\frac{R^2(s)}{b(s)}+1\right)}}\,e^{-\frac{\left(x+\frac{c(s)}{b(s)}\right)^2}{-\frac{R^2(s)}{b(s)}+1}},
\end{eqnarray*}
\begin{eqnarray*}
 Q_-^*(x,s)=\frac{1}{\sqrt{{\pi}\left(-\frac{R^2(s)}{b(s)}+1\right)}}\,e^{-\frac{\left(x+\frac{c(s)}{b(s)}\right)^2}{-\frac{R^2(s)}{b(s)}}} \le  Q^*(x,s)  \le  Q_+^*(x,s)=\frac{1}{\sqrt{-\frac{\pi R^2(s)}{b(s)}}}\,e^{-\frac{\left(x+\frac{c(s)}{b(s)}\right)^2}{-\frac{R^2(s)}{b(s)}+1}},
\end{eqnarray*}
and therefore
\begin{eqnarray*}
&&B_0^*(x,s) - \frac{A_{10}}{A_{11}} Q_-^*(x,s) \le    p^*(x,s)\le  B_0^*(x,s) - \frac{A_{10}}{A_{11}} Q_+^*(x,s),\quad A_{10}\ge 0,\\
&&B_0^*(x,s) - \frac{A_{10}}{A_{11}} Q_+^*(x,s) \le    p^*(x,s)\le  B_0^*(x,s) - \frac{A_{10}}{A_{11}} Q_-^*(x,s),\quad A_{10}< 0.
 \end{eqnarray*}
 Note that the exact expression for $Q^*(x,s)$ is given through the integral of the Whittaker function.

\medskip

c.
Let us choose the initial data such that the dynamics in the case $b(s)<0$ can be found explicitly. Namely,
\begin{eqnarray}\label{IDDel}
 \Phi(x,s)=\delta (x).
\end{eqnarray}
 then $\hat g_0(\mu)= 1$,  $\hat g_1(\mu)= 0$. From \eqref{9b} we find $ a_0(t,\mu)=1$, $
  a_1(t,\mu)=\frac{A_{10}}{A_{11}}(e^{A_{11}t}-1)$.

Further, from \eqref{11b} and  \eqref{16b} we have
\begin{eqnarray*}\label{17e}
&& p(t,x,s)=\\&&\frac{1}{\sqrt{-\frac{\pi R^2(s)}{b(s)}(1-e^{2b(s)t})}}\, \left(1+ \frac{A_{10}}{A_{11}}\, \left(e^{A_{11} t}-1\right)\,X(s)\right)\,e^{-\frac{\left(x+\frac{c(s)}{b(s)}(1-e^{b(s)t})\right)^2}{-\frac{R^2(s)}{b(s)}(1-e^{2b(s)t})}}.\nonumber
\end{eqnarray*}

\medskip

\subsubsection{Stepwise Gaussian initial distribution}\label{3.3.2}

 Let the initial data be
\begin{eqnarray}\label{ID1}
 \Phi(x,s)=\frac{1} {\sqrt{\pi}}\,e^{-(x-m(s))^2},\quad m(s)= \left\{\begin{array}{cc}m_1,\, s\in [0,\frac12]\\
 m_2, \, s \in (\frac12,1].
  \end{array}\right.
\end{eqnarray}
We can expand \eqref{ID1} the into  series \eqref{14}. To use the results of Sec.\ref{Stwomodes} we consider two first member of the expansion only:
\begin{eqnarray*}
&& \bar \Phi(x,s)=g_0(x)+ g_1(x) X_1(s).\label{Phibar}
\end{eqnarray*}
If $X_1(s)=\cos \pi s$, then
\begin{eqnarray}
&& g_0(x)=\frac{1} {2\sqrt{\pi}}\,  (e^{-(x-m_1)^2}+e^{-(x-m_2)^2}),\qquad
g_1(x)= \frac{2} {\pi^\frac32}\, (e^{-(x-m_1)^2}-e^{-(x-m_2)^2}).\nonumber
\end{eqnarray}
Evidently, $ \bar \Phi(x,s)=g_0(x)+ g_1(x) \cos \pi s \ne \Phi(x,s) $.
In this case
\begin{eqnarray*}
&&  B_0(t,x,s)= \frac{1}{2\sqrt{\pi(1+2 t R^2(s)})} (e^{-\frac{(x-m_1-c(s) t)^2}{1+2 t R^2(s)}}+e^{-\frac{(x-m_2-c(s) t)^2}{1+2 t R^2(s)}}), \\
&&B_1(t,x,s)= \frac{2}{\sqrt{\pi^3(1+2 t R^2(s))} } (e^{-\frac{(x-m_1-c(s) t)^2}{1+2 t R^2(s)}}-e^{-\frac{(x-m_2-c(s) t)^2}{1+2 t R^2(s)}}),
\end{eqnarray*}

\medskip

If we choose $X_1(s)=H_1(s)$, then
\begin{eqnarray}
&& g_0(x)=\frac{1} {2\sqrt{\pi}}\,  (e^{-(x-m_1)^2}+e^{-(x-m_2)^2}),\qquad
g_1(x)= \frac{1} {2\sqrt{\pi}}\, (e^{-(x-m_1)^2}-e^{-(x-m_2)^2}),\nonumber
\end{eqnarray}
end $ \bar \Phi(x,s)=g_0(x)+ g_1(x) H_1(s) = \Phi(x,s) $.
In this case
\begin{eqnarray*}
&&  B_0(t,x,s)= \frac{1}{2\sqrt{\pi(1+2 t R^2(s)})} (e^{-\frac{(x-m_1-c(s) t)^2}{1+2 t R^2(s)}}+e^{-\frac{(x-m_2-c(s) t)^2}{1+2 t R^2(s)}}), \\
&&B_1(t,x,s)= \frac{1}{2\sqrt{\pi (1+2 t R^2(s))} } (e^{-\frac{(x-m_1-c(s) t)^2}{1+2 t R^2(s)}}-e^{-\frac{(x-m_2-c(s) t)^2}{1+2 t R^2(s)}}),
\end{eqnarray*}

\medskip

 a. For $b=0$ from  \eqref{10}  and  \eqref{16} we have the two-modes approximation of density $p(t,x,s)$ as
\begin{eqnarray*}\label{p1}
&& \bar p(t,x,s)=\\
&&B_0(t,x,s)+\left(e^{A_{11} t} B_1(t,x,s)+\frac{A_{10}}{A_{11}} \left(e^{A_{11} t}-1\right) B_0(t,x,s)\right)\, X_1(s).\nonumber
 \end{eqnarray*}
We see that in the case $A_{11}<0$, the limit behavior do not depend on $B_1$,
\begin{eqnarray}\label{p11}
&& \bar p(t,x,s)\sim\left(1-\frac{A_{10}}{A_{11}} \right) \, B_0(t,x,s)\, X_1(s),\quad t\to \infty.\nonumber
 \end{eqnarray}

 \medskip

 b. For $b\ne 0$ analogically to the case b, Sec.\ref{S321}, we get from \eqref{11b} and  \eqref{16b}
 \begin{eqnarray*}\label{17B00}
&&B_0(t,x,s)=
\frac{e^{-\frac{\left(x-m_1+(1-e^{b(s)t})\frac{c(s)}{b(s)}\right)^2}{\frac{R^2(s)}{b(s)}( e^{2 b(s) t-1})+1}}+e^{-\frac{\left(x-m_2+(1-e^{b(s)t})\frac{c(s)}{b(s)}\right)^2}{\frac{R^2(s)}{b(s)}( e^{2 b(s) t-1})+1}}}{2 \sqrt{\pi\left(\frac{R^2(s)}{b(s)}(e^{2 b(s) t}-1)+1\right)}},
\end{eqnarray*}
 and
\begin{eqnarray*}
&& B_1(t,x,s) = \frac{A_{10}}{A_{11}} (e^{A_{11} t}-1)  Q_b(t,x,s),
\end{eqnarray*}
where $Q_b(t,x,s)$ can be estimated from both sides as in Sec.\ref{S321}, we do not repeat these rather cumbersome estimates.

\medskip

 It can be readily shown that if $b<0$, $A_{11}<0$,   then the limiting distribution of the density function contains two peaks, and structure of this distribution is actually determined only by the coefficient $g_0(x)$.

\medskip

c. To find explicit  solution in the case $b(s)<0$ we choose the initial data
\begin{eqnarray}\label{IDel}
 \Phi(x,s)=\delta (x-m(s)),
\end{eqnarray}
where $m(s)$ is given in \eqref{ID1}.

If  $X_1(s)=H_1(s)$, then
\begin{eqnarray}
&& g_0(x)=\frac{1} {2}\,  (\delta (x-m_1)+\delta (x-m_2)),\qquad
g_1(x)= \frac{1} {2}\,  (\delta (x-m_1)-\delta (x-m_2)).\nonumber
\end{eqnarray}

Further, from \eqref{11b} and  \eqref{16b} we have
\begin{eqnarray}\label{pp1}
&&  p(t,x,s)=\\
&&B_0(t,x,s)+\left(e^{A_{11} t} B_1(t,x,s)+\frac{A_{10}}{A_{11}} \left(e^{A_{11} t}-1\right) B_0(t,x,s)\right)\, X_1(s),\nonumber
 \end{eqnarray}
\begin{eqnarray*}
&&B_0(t,x,s)=\frac{e^{-\frac{\left(x-m_1+\frac{c(s)}{b(s)}(1-e^{b(s)t})\right)^2}{-\frac{R^2(s)}{b(s)}(1-e^{2b(s)t})}}+
e^{-\frac{\left(x-m_2+\frac{c(s)}{b(s)}(1-e^{b(s)t})\right)^2}{-\frac{R^2(s)}{b(s)}(1-e^{2b(s)t})}} }{2\,\sqrt{-\frac{\pi R^2(s)}{b(s)}(1-e^{2b(s)t})}},\\
&&B_1(t,x,s)=\frac{e^{-\frac{\left(x-m_1+\frac{c(s)}{b(s)}(1-e^{b(s)t})\right)^2}{-\frac{R^2(s)}{b(s)}(1-e^{2b(s)t})}}-
e^{-\frac{\left(x-m_2+\frac{c(s)}{b(s)}(1-e^{b(s)t})\right)^2}{-\frac{R^2(s)}{b(s)}(1-e^{2b(s)t})}} }{2\,\sqrt{-\frac{\pi R^2(s)}{b(s)}(1-e^{2b(s)t})}}.
\end{eqnarray*}

\medskip

Pic.1 and 2 present the behavior of densities for the main and hidden states for the cases $b=0$  for initial data \eqref{ID1} and $b<0$ for initial data \eqref{IDel},
decomposition is over the Haar basis $H_n(s)$, $b(s)= \left\{\begin{array}{cc}b_1,\, s\in [0,\frac12]\\
 b_2, \, s \in (\frac12,1]\end{array}\right.$, $R(s)= \left\{\begin{array}{cc}R_1,\, s\in [0,\frac12]\\
 R_2, \, s \in (\frac12,1]\end{array}\right.$, the parameters are taken as
\begin{eqnarray}\label{param}
&&\lambda_1 = 1, \, \lambda_2 = 2, \, R_1 = 1, \,  R_2 = 2, \\
&&  m_1 = 5, m_2 = -5,\,b_1 = -0.5,\, b_2 = -1, \,   c =1.\nonumber
\end{eqnarray}

\begin{center}
\begin{figure}[htb!]
\hspace{-1cm}
\hspace{1cm}
\begin{minipage}{0.4\columnwidth}
\includegraphics[scale=0.35]{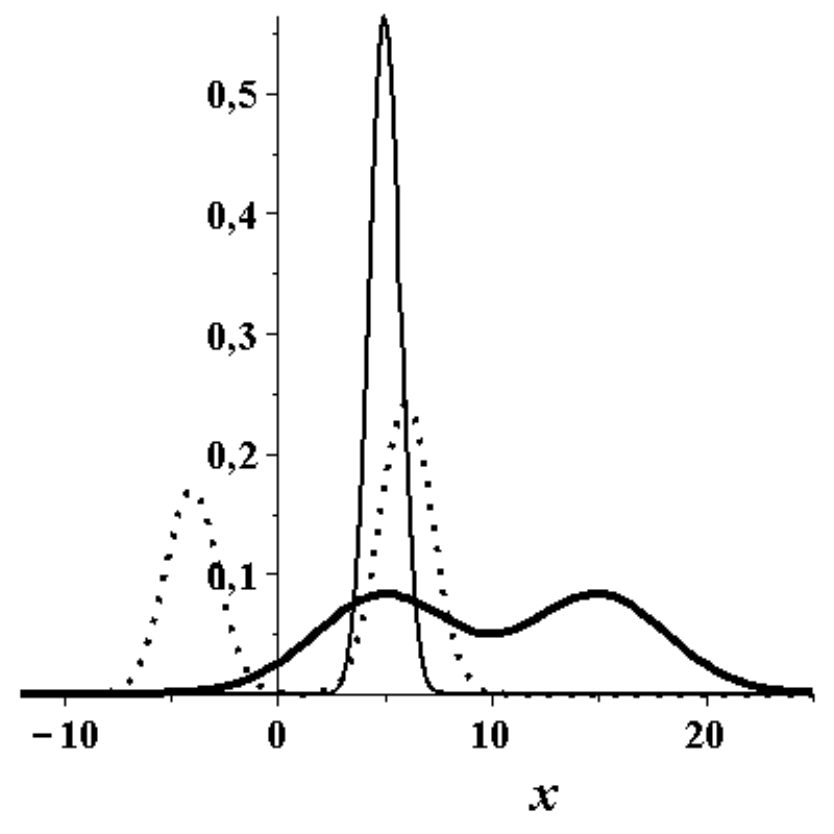}
\end{minipage}
\hspace{1.5cm}
\begin{minipage}{0.4\columnwidth}
\includegraphics[scale=0.35]{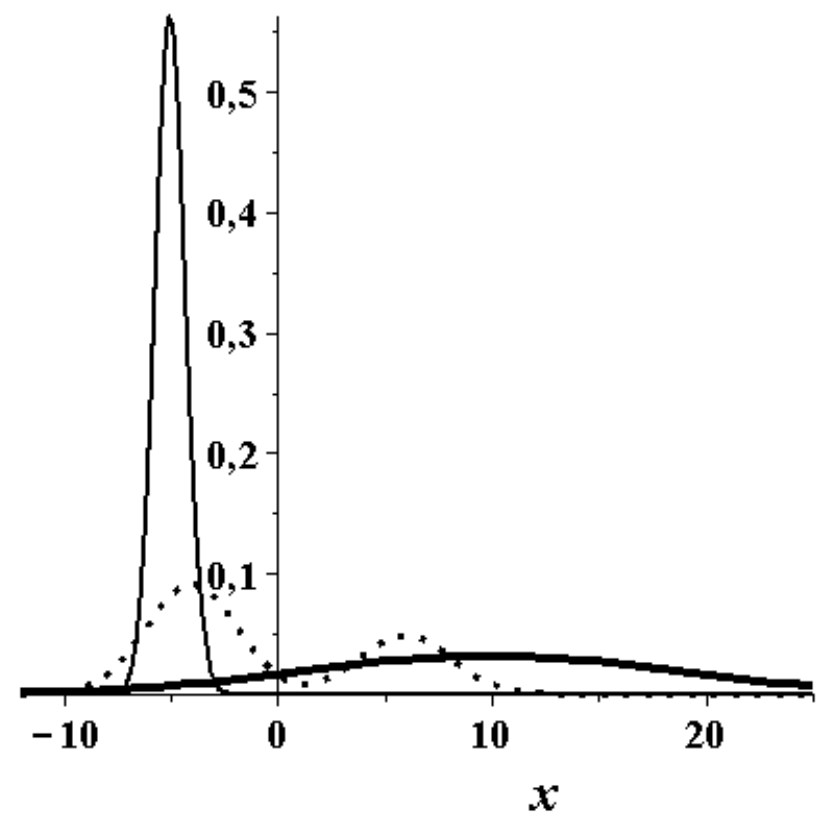}
\end{minipage}
\caption{The dynamics of densities of main (left) and hidden (right) states for $b=0$, other parameters are given in \eqref{param};
$t=0$ (thin solid line), $t=1$ (dots), $t=10$ (thick solid line). The initial data are \eqref{ID1}.
}\label{Pic1}
\end{figure}
\end{center}

\begin{center}
\begin{figure}[htb!]
\hspace{-1cm}
\hspace{1cm}
\begin{minipage}{0.4\columnwidth}
\includegraphics[scale=0.35]{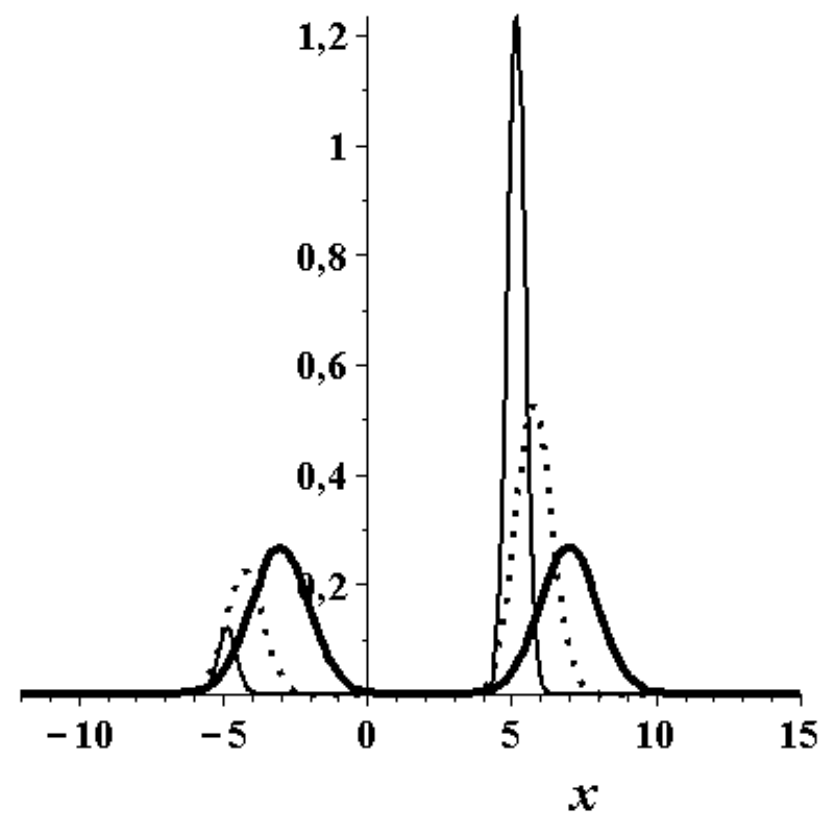}
\end{minipage}
\hspace{1.5cm}
\begin{minipage}{0.4\columnwidth}
\includegraphics[scale=0.35]{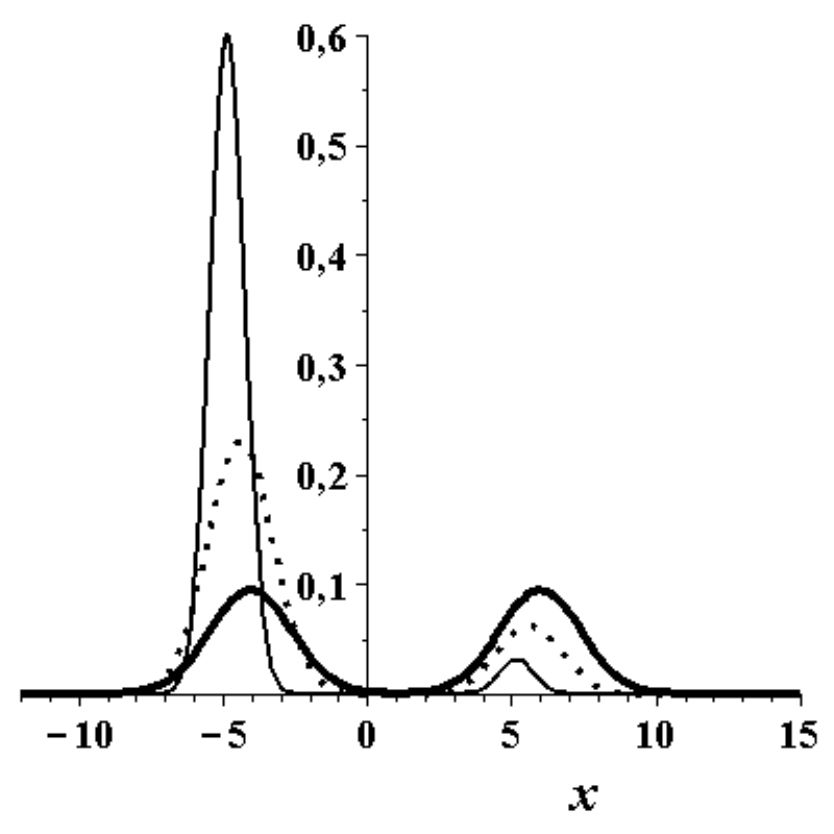}
\end{minipage}
\caption{The dynamics of densities of main (left) and hidden (right) states for $b<0$, the parameters are given in \eqref{param}; $t=0.1$ (thin solid line), $t=0.5$ (dots), $t=100$ (thick solid line). The initial data are  \eqref{IDel}.
}\label{Pic1}
\end{figure}
\end{center}

As can be seen from \eqref{9b}, the time-independent solution of  \eqref{3} depends on $g_0(x)$, that is, on the initial conditions \eqref{4}. This is an unusual  result. In particular, in the case of the Gaussian initial data \eqref{IDGauss}, this equilibrium position cannot be found explicitly; it is expressed through the Whittaker function. For the delta-function initial data \eqref{IDDel}, this solution is found explicitly; it can also be obtained by passing to the limit as $t\to\infty$ in \eqref{pp1}.

\section{A four-state switching diffusion, $M=3$}

To apply our result to higher state switching diffusion based in the  continuous-time multifractal model, see, e.g.\cite{Lux}.
 On the level of hierarchy equal to $k$, $k\in \mathbb N$, the states can be enumerated by groups of $k$ elements consisting of $0$ and $1$;  it allows for a total of ${2^k}$ Markov states. For $k=1$, the states belong to the set $s=\{0,1\}$, the coefficients of diffusion are $m_0$ and $m_1=2-m_0$, $0<m_0<2$. Thus, it is the situation considered in Sec.\ref{2st}.
For $k = 2$ the states can be enumerated as $ s= \{00, 10, 01, 11\}$.
Since the hierarchical
components enter multiplicatively,  the state dependent variances are
$m_0 m_0, (2 - m_0) m_0, m_0 (2 - m_0),  (2 - m_0)(2 - m_0)$.
The system of four partial differential equations governing the  densities $u_{ml}(t,x)$, $m,l\in {0,1}$ is
\begin{eqnarray}\label{uuuu}
&&(u_{00})_t=\frac12 r_0^2\,(u_{00})_{xx}- (\lambda_1+\lambda_2) u_{00}+\lambda_1 u_{10}+\lambda_2 u_{01},\nonumber\\
&&(u_{10})_t=\frac12 r_0 (2-r_0)\,(u_{10})_{xx}+ \lambda_1 u_{00} - (\lambda_1+\lambda_2) u_{10} + \lambda_2 u_{11},\\
&&(u_{01})_t=\frac12 r_0 (2-r_0)\,(u_{01})_{xx}+ \lambda_2 u_{00} - (\lambda_1+\lambda_2) u_{01} + \lambda_1 u_{11},\nonumber\\
&&(u_{11})_t=\frac12 (2-r_0)^2\,(u_{11})_{xx}+ \lambda_2 u_{10} + \lambda_1 u_{01} - (\lambda_1+\lambda_2) u_{11},\nonumber
\end{eqnarray}
with constant $\lambda_i$, $i=1, 2$, $0<\lambda_2<\lambda_1$.

Let us show that this system can be explicitly solved in our framework. Indeed, for the convenience we enumerate states as
$w(t,0)=u_{00}$, $w(t,1)=u_{10}$, $w(t,2)=u_{01}$, $w(t,3)=u_{11}$, and $u_{00}$ corresponds to the main state,
\begin{equation}\label{qji4}
q_{ji}=\left(\begin{array}{cccc}-(\lambda_1+\lambda_2)
 & \lambda_1 & \lambda_2& 0\\  \lambda_2&-(\lambda_1+\lambda_2)& 0 & \lambda_1\\
 \lambda_1 & 0 &-(\lambda_1+\lambda_2)& \lambda_2\\
 0&  \lambda_2& \lambda_1&-(\lambda_1+\lambda_2)
  \end{array}\right),
\end{equation}
$i,j=0,1,2,3,$
\begin{eqnarray*}\label{RRR}
  R^2(s)= \left\{\begin{array}{cc}r_0^2,\, s\in [0,\frac14]\\
 r_0 (2-r_0), \, s \in (\frac14,\frac12]\\ r_0 (2-r_0), \, s \in (\frac12,\frac34] \\ (2-r_0)^2, \, s \in (\frac34,1].
  \end{array}\right.
\end{eqnarray*}

It is convenient to use as a basis $X_n(s)$, $n=0,1,\dots$, the Haar wavelet basis $H_n(s)$ as in Sec.\ref{Xex}. We need only 4 first function of the orthonormal basis, that is
\begin{eqnarray*}\label{X4}
&&H_0(s)=1, \quad H_1(s)=\phi(s), \quad H_2(s)=\sqrt{2}\, \phi \left(2^{\frac{1}{2}} s \right), \\
  &&H_3(s)=\sqrt{2}\, \phi \left(2^{\frac{1}{2}} s-1 \right).\nonumber
\end{eqnarray*}
It can be readily found that the matrix of Fourier coefficients (see \eqref{7}) is
\begin{equation}\label{Anm4}
A_{nm}=\left(\begin{array}{cccc} 0 & 0& 0&0\\
0& -\frac14 (\lambda_1+\lambda_2)&-\frac{\sqrt{2}}{8} (\lambda_1-\lambda_2)&-\frac{\sqrt{2}}{8} (\lambda_1-\lambda_2)\\
0& \frac{\sqrt{2}}{8} (\lambda_1-\lambda_2)&   -\frac38 (\lambda_1+\lambda_2)& \frac18 (\lambda_1+\lambda_2)\\
0& \frac{1}{2} (\lambda_1-\lambda_2)&\frac18 (\lambda_1+\lambda_2)&   -\frac38 (\lambda_1+\lambda_2)
  \end{array}\right),
\end{equation}
 $n,m=0,1,2,3$.
We use this matrix to find coefficients $a_k(t,\mu)$, $k=0,1,2,3,$ from system \eqref{9} (or \eqref{9b})

We choose the initial data as
\begin{eqnarray}\label{ID4}
 \Phi(x,s)=\frac{1} {\sqrt{\pi}}\,e^{-(x-m(s))^2},\quad m(s)= \left\{\begin{array}{cc}m_0,\, s\in [0,\frac14]\\
 m_1, \, s \in (\frac14,\frac12]\\ m_2, \, s \in (\frac12,\frac34] \\ m_3, \, s \in (\frac34,1].
  \end{array}\right.
\end{eqnarray}
Then
\begin{eqnarray*}
&&\Phi(x,s)=g_0(x)+ g_1(x) H_1(s) + g_2(x) H_2(s)+ g_3(s) H_3,\\
&& g_0(x)=\frac{1} {4\sqrt{\pi}}\,  \left(e^{-(x-m_0)^2}+e^{-(x-m_1)^2}+e^{-(x-m_2)^2}+e^{-(x-m_3)^2}\right),\\
&& g_1(x)=\frac{1} {4\sqrt{\pi}}\,  \left(e^{-(x-m_0)^2}-e^{-(x-m_1)^2}+e^{-(x-m_2)^2}-e^{-(x-m_3)^2}\right),\\
&& g_2(x)=\frac{1} {\sqrt{2\pi}}\,  \left(e^{-(x-m_0)^2}-e^{-(x-m_1)^2}\right),\qquad g_3(x)=\frac{1} {\sqrt{2\pi}}\,  \left(e^{-(x-m_2)^2}-e^{-(x-m_3)^2}\right).
\end{eqnarray*}

Since $b=c=0$, we have to solve \eqref{9} subject to
$a_{k}(0,\mu)=\hat g_k(\mu)$.
The eigenvalues of $A_{nm}$, $n,m=1, 2, 3$, a minor of \eqref{Anm4} that determines the solution,
 are $-\frac12(\lambda_1+\lambda_2)$, $-\frac14(\lambda_1+\lambda_2)(1\pm i)$.
 In can be computed that
\begin{eqnarray*}
&&  B_0(t,x,s)= \frac{1}{4\sqrt{\pi(1+2 t R^2(s)})} \left(e^{-\frac{(x-m_0)^2}{1+2 t R^2(s)}}+e^{-\frac{(x-m_1)^2}{1+2 t R^2(s)}}+e^{-\frac{(x-m_2)^2}{1+2 t R^2(s)}}+e^{-\frac{(x-m_3)^2}{1+2 t R^2(s)}}\right), \\
&&B_1(t,x,s)=\frac{\sqrt{2}}{2} e^{-\frac{\lambda_1+\lambda_2}{4} t } \,\left(\sqrt{2}\,\alpha_1 \cos \left(\frac{\lambda_2-\lambda_1}{4}
 t\right) +(\alpha_2+\alpha_3)\,\sin \left(\frac{\lambda_2-\lambda_1}{4}t\right) \right),\\
 &&B_2(t,x,s)=\frac12 (\alpha_2-\alpha_3) e^{-\frac{\lambda_1+\lambda_2}{2} t}
  -\frac12 \alpha_1  e^{-\frac{\lambda_1+\lambda_2}{4} t }\sin \left(\frac{\lambda_2-\lambda_1}{4} t\right)+ \\
  &&  \phantom{B_2(t,x,s)=}  \frac12 (\alpha_2-\alpha_3)\,e^{-\frac{\lambda_1+\lambda_2}{4} t }\cos \left(\frac{\lambda_2-\lambda_1}{4} t\right)\\
  &&B_3(t,x,s)=-\frac12 (\alpha_2-\alpha_3) e^{-\frac{\lambda_1+\lambda_2}{2} t}
  -\frac12 \alpha_1  e^{-\frac{\lambda_1+\lambda_2}{4} t }\sin \left(\frac{\lambda_2-\lambda_1}{4} t\right)+ \\
  &&  \phantom{B_2(t,x,s)=}  \frac12 (\alpha_2-\alpha_3)\,e^{-\frac{\lambda_1+\lambda_2}{4} t }\cos \left(\frac{\lambda_2-\lambda_1}{4} t\right),
\end{eqnarray*}
where
\begin{eqnarray*}
&&  \alpha_1(t,s)=\frac{1}{4\sqrt{\pi(1+2 t R^2(s)})} (e^{-\frac{(x-m_0)^2}{1+2 t R^2(s)}}+e^{-\frac{(x-m_1)^2}{1+2 t R^2(s)}}-e^{-\frac{(x-m_2)^2}{1+2 t R^2(s)}}-e^{-\frac{(x-m_3)^2}{1+2 t R^2(s)}}),\\
&&  \alpha_2(t,s)=\frac{1}{2 \sqrt{2\pi(1+2 t R^2(s)})} (e^{-\frac{(x-m_0)^2}{1+2 t R^2(s)}}-e^{-\frac{(x-m_1)^2}{1+2 t R^2(s)}}),\\
 &&  \alpha_3(t,s)=\frac{1}{2 \sqrt{2\pi(1+2 t R^2(s)})} (e^{-\frac{(x-m_2)^2}{1+2 t R^2(s)}}-e^{-\frac{(x-m_3)^2}{1+2 t R^2(s)}}).
\end{eqnarray*}
Thus,
\begin{eqnarray*}\label{p4}
&&p(t,x,s)=B_0(t,x,s)+B_1(t,x,s) H_1(s)+B_2(t,x,s) H_2(s)+B_3(t,x,s) H_3(s).
\end{eqnarray*}

Pic.3 present the behavior of densities for the main and 3 hidden states for the cases for initial data \eqref{ID4},
decomposition is over the Haar basis $H_n(s)$,   the parameters are taken as
\begin{eqnarray}\label{param4}
&&\lambda_1 = 0.4, \, \lambda_2 = 0.2, \, r_0 = 1.7, \,
 m_0=-10,  m_2 = 1, m_2 = -5, m_3=10.
\end{eqnarray}
Recall that all states  tend to zero as $t\to\infty$.

\begin{center}
\begin{figure}[htb!]
\hspace{-1cm}
\hspace{1cm}
\begin{minipage}{0.4\columnwidth}
\includegraphics[scale=0.35]{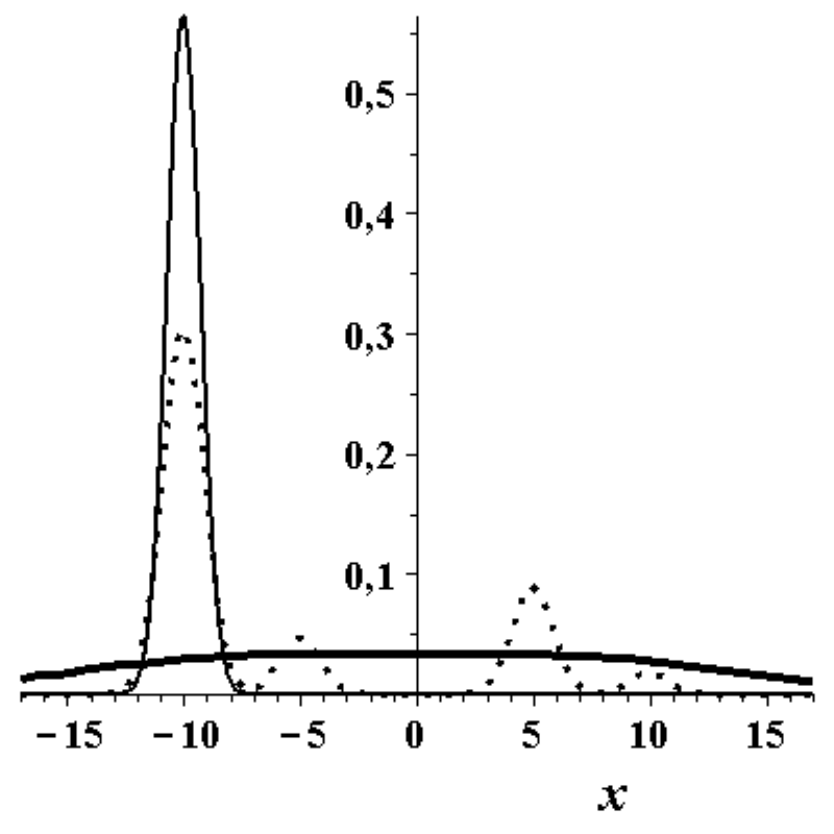}
\end{minipage}
\hspace{1.5cm}
\begin{minipage}{0.4\columnwidth}
\includegraphics[scale=0.35]{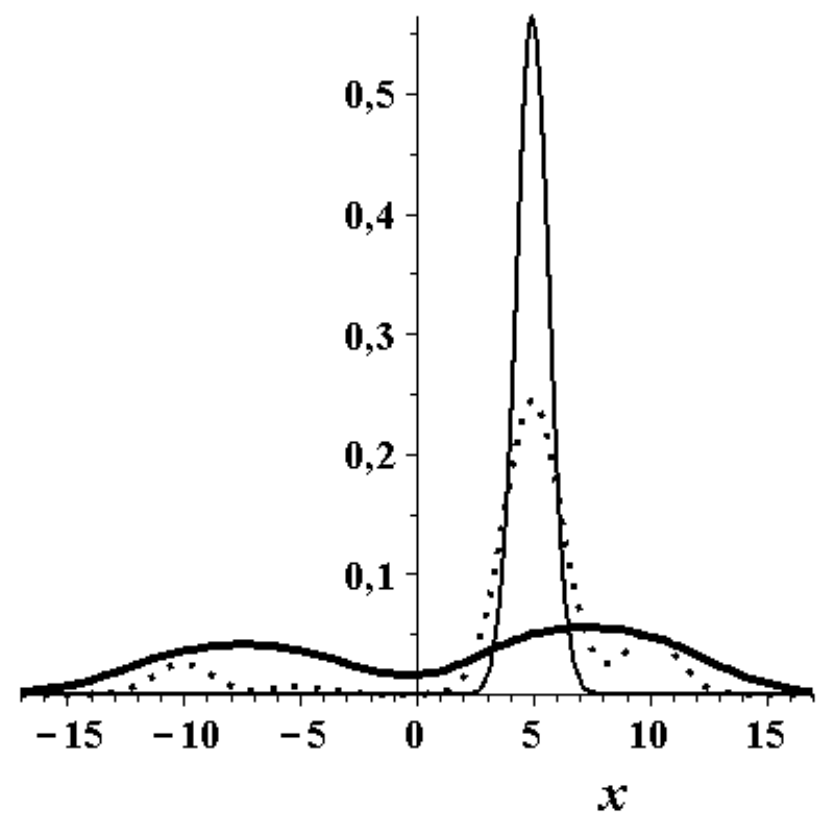}
\end{minipage}
\hspace{1.5cm}
\begin{minipage}{0.4\columnwidth}
\includegraphics[scale=0.35]{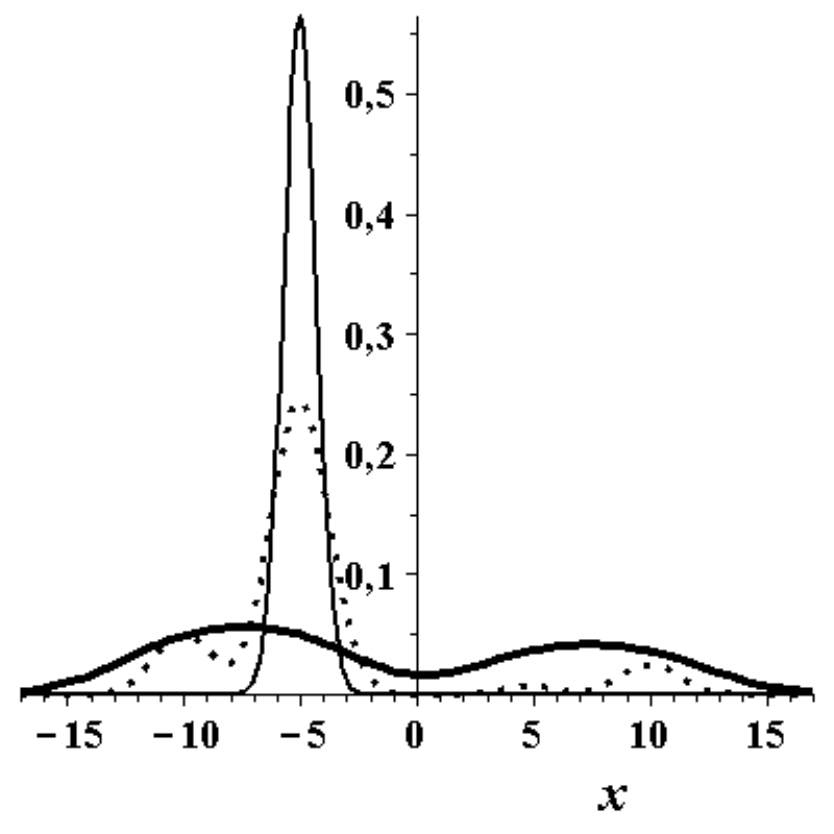}
\end{minipage}
\hspace{1.5cm}
\begin{minipage}{0.4\columnwidth}
\includegraphics[scale=0.35]{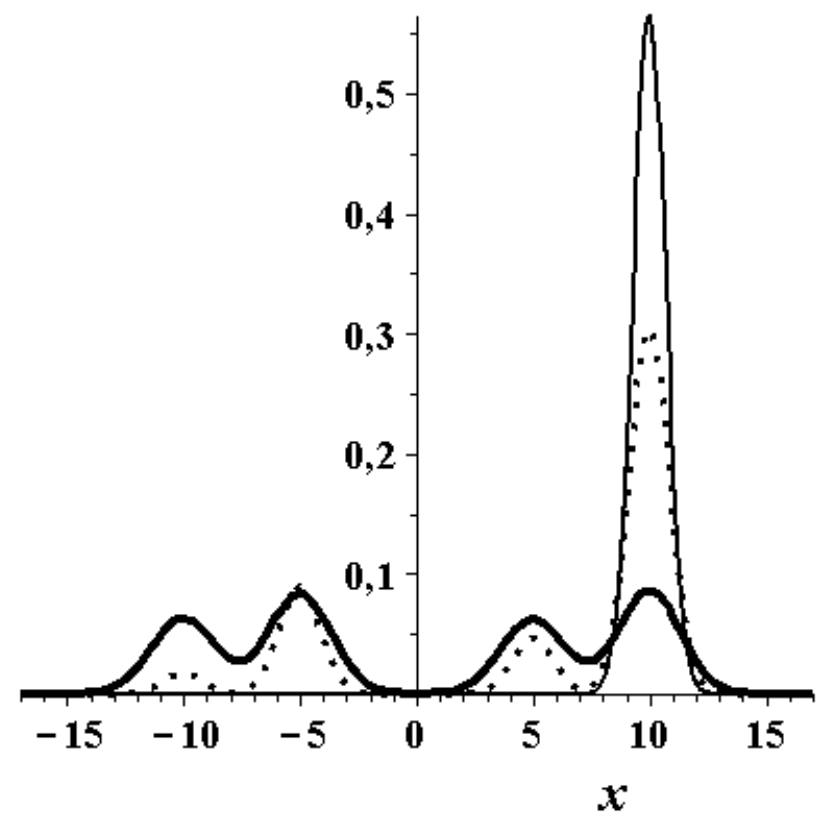}
\end{minipage}
\caption{The dynamics of densities of main (top left) and hidden states 1, 2, 3 (top right, bottom left, bottom right, respectively) for parameters given in \eqref{param4};
$t=0$ (thin solid line), $t=3$ (dots), $t=15$ (thick solid line).
}\label{Pic1}
\end{figure}
\end{center}

\begin{rem} \rm{
1. Note that $p(t,x,s)=B_0(t,x,s) + \gamma(t)$, $\gamma(t)\to 0$ as  $t\to\infty$, but the character of $\gamma(t)$ is  oscillatory.

2. The structure if the matrix \eqref{qji4} is quite simple (it contains only 2 free parameters), which  determines the relative simplicity of solution \eqref{p4}.
For an arbitrary matrix $q_{ij}$, having $q$ - property, having 12 free parameters,   the eigenvalues of $A_{nm}$, $n,m=1,2,3$, also have  a negative real part. Indeed, the eigenvalues are roots of equation
\begin{eqnarray*}
\xi^3 + k_2 \xi^2 + k_1 \xi + k_0=0,
\end{eqnarray*}
where the coefficients $k_i$ depend on $q_{ij}$ in a very complex way (the result obtained using a computer algebra package).
However, all coefficients are positive, this  implies that the real part of all roots are negative. A similar result should be true for any dimension, this is a consequence of the maximum principle for equation \eqref{3}.

3. If we add an advective term $\beta$ to the model \eqref{uuuu}, we obtain for $b<0$ a nonzero structure  as $t\to\infty$, similar to what we saw in the case of $M=1$, Sec.\ref{3.3.2}. For $M=3$ it  consists of four spikes.

4. At hierarchy levels $k=3, 4,$ etc., the dimension of matrix $q_{ij}$ is very high, $(2^k \times 2^k)$, and it seems impossible to conduct any analytical reasoning.
 However, this matrix is very sparse,
 that is, it contains a large number of zeros,
 and becomes increasingly sparse
with increasing $k$.
Furthermore, it has a convenient block structure (see its description in \cite{Lux}, Sec.4.2). At hierarchy level $k$, there are only $k$ free parameters, so results can still be obtained.}
\end{rem}

\section{Discussion}
We study an integro-differential equation, which is a generalization of the Fokker-Planck system of equations for each of the states between which switching occurs, to the case of a continuum of states. We show that this approach can yield exact solutions to the equation in the form of a series or a finite sum. In particular, with a certain choice of kernels, solving the integro-differential equation can yield a solution for a discrete model. The key is the choice of an orthonormal basis over which to expand the kernel and initial data. For the transition to a discrete model, it is convenient to use the Haar basis. Using this same methodology, one can construct a solution not only to the Fokker-Planck equation for diffusion with switching, but also to the generalized Ornstein-Uhlenbeck model. For the latter model, the solution converges to a nontrivial constant state, which can also be constructed explicitly.
This is a structure consisting of several peaks.  Similar structures appeared as given stationary distributions \cite{Tretyakov}, but with a special selection of the switching matrix $q_{ji}$. The stability of regime-switching diffusions was studied in  \cite {Khasminskii}.
Note that numerical methods for solving stochastic differential equations are successfully used to study the properties of the forward and backward Kolmogorov equations (e.g., \cite{Mao_Yuan}, \cite{Bao} \cite{Tretyakov} and references therein).

We also note that our method can be easily generalized to the multivariate case in $x$. However, obtaining an explicit solution is limited by the type of initial data: it must be a delta-function or a Gaussian distribution in $x$. Incidentally, such distributions are the most common in applications.

\section*{Acknowledgments}
A.B. (conceptualisation) was supported by Moscow Center for Fundamental and Applied Mathematics.
O.R. (writing, investigation, pictures) was supported by RSF grant 23-11-00056 through RUDN University.

\end{document}